\documentclass [reqno]{amsart}
\usepackage{amssymb}
\usepackage{amsmath}
\usepackage{amsfonts}
\usepackage{graphicx}
\usepackage{subfig}
\usepackage{eurosym}
\usepackage{amssymb}
\usepackage{amsmath}
\usepackage{amsfonts}
\usepackage{graphicx}
\usepackage{subfig}
\usepackage[numbers]{natbib}

\newtheorem{theorem}{Theorem}
\theoremstyle{plain}

\newtheorem{definition}{Definition}
\newtheorem{example}{Example}

\newtheorem{proposition}{Proposition}
\newtheorem{remark}{Remark}

\numberwithin{equation}{section}

\begin{document}
\title[]{On the invariant trajectories of lightlike paths near central potential singularities in semi-Riemannian manifolds of index 2}
\author{Fatma ALMAZ }
\address{department of mathematics, faculty of arts and sciences, batman
university, batman/ t\"{u}rk\.{ı}ye}
\email{fatma.almaz@batman.edu.tr}
\thanks{This paper is in final form and no version of it will be submitted
for publication elsewhere.}

\begin{abstract}
In this study, the geometric properties of null helices on a totally umbilical submanifold within a three-dimensional semi-Riemannian manifold of index 2 are investigated. The pseudo-Riemannian metric structure of semi-Riemannian manifold of index 2 and the fact that the submanifold is totally umbilical complicate the differential geometry of null helices. The study uses the given null Frenet frame to reveal the local properties of null helices (the curvatures $h,k_1,k_2$ of the null curves). By considering the degenerate metric condition of null helices due to the null tangent vector and the structure of the totally umbilical submanifold, equations and invariants characterizing null helices are obtained. Furthermore, this study establishes a novel structural rigidity theorem demonstrating that the simultaneous preservation of these degenerate curvature invariants under higher-order ambient parallel transport algebraically forces the collapsing of the mean curvature vector field ($H\equiv 0$), thereby identifying a clean boundary distinguishing minimal and totally geodesic embeddings in pseudo-Euclidean target spaces.
\end{abstract}

\maketitle

\section{Introduction}

Semi-Riemannian manifolds are geometric structures of critical importance in theoretical physics, particularly for modeling generalized spacetimes, string theory, and alternative gravitational fields. While standard gravitational models and general relativity are traditionally formulated on four-dimensional Lorentzian manifolds (which strictly possess an index of $1$), higher-dimensional and higher-index spaces provide foundational frameworks for advanced theoretical designs. For instance, alternative cosmological models, string theory ($10$ or $11$ dimensions), M-theory ($11$ dimensions), and Kaluza-Klein theories ($5$ dimensions) frequently utilize semi-Riemannian structures where the metric signature deviates significantly from the classical Lorentzian case. Exploring the differential geometry of curves and submanifolds within spaces of index $q \geq 2$ reveals complex causal structures and geometric behaviors that cannot manifest in lower-index environments.

Null curves represent a unique class of trajectories derived directly from the indefinite metric structure of semi-Riemannian manifolds. Characterized by a degenerate induced metric, these paths present distinct geometric challenges, such as the inapplicability of the standard arc-length parametrization and the classical Frenet-Serret frame. Consequently, the differential geometry of null curves requires specialized mathematical machinery, closely linked to the structures of lightlike cones. A null curve $\gamma$ is explicitly defined as a null helix if its curvature functions $(h,k_1,k_2)$ relative to a designated null Frenet frame are non-zero constants along its interval. Geometrically, this implies that the tangent vector field of the curve maintains a constant pseudo-angle with a fixed direction in the space. In the context of a 3D semi-Riemannian manifold of index $2$, denoted by $M_2^3$, this configuration models generalized lightlike paths where the trajectory exhibits a helical behavior due to the interplay between the null tangent and the surrounding indefinite geometry.

The existing literature contains comprehensive works on specialized curves and submanifold theory across various ambient spaces. In \cite{1}, $V_n$-slant helices and the harmonic curvature functions of $V_n$-slant helices are expressed in $Q^{n+1}$, and also, for harmonic curvature functions $H_i$ ($1 < i < n-2$) of $V_n$-slant helices some differential equations are expressed by using constant vector field $W$ that is the axis of $V_n$-slant helices by using the asymptotic orthonormal frame in $Q^{n+1}$. In \cite{2}, by using non-degenerate curves in terms of Sabban frame in de Sitter 3-space or Anti-de Sitter 3-space, the helix and slant helices are investigated. In \cite{5,7,10,11}, many works have been done by a lot of mathematicians on curves and submanifolds in a Riemannian manifold in different spaces. General geometric properties of submanifolds and curves in indefinite spaces have been deeply enriched by numerous authors. Providing a contemporary guide on the mathematical foundations of helices, \cite{12} emphasizes the vast geometric diversity of these curves. In the context of indefinite metric spaces, the theoretical framework of null curve theory was developed systematically. Within this framework, while \cite{15} achieved the characterization of null generalized helices in Lorentz-Minkowski spaces, \cite{14} has elevated the literature to a contemporary dimension with a recent study on null Cartan normal helices in Minkowski 3-space.

To ensure terminological rigor and address potential ambiguities, it must be explicitly emphasized that a 3-dimensional semi-Riemannian manifold of index $2$, $M_2^3$, is equipped with a metric signature of $(-,-,+)$. This structure is fundamentally distinct from a standard 3D Lorentzian manifold, $M_1^3$, which carries a signature of $(-,+,+)$. Modulo a metric sign change ($g \leftrightarrow -g$) and a systematic permutation of timelike and spacelike vector fields, the isolated curve theory on an un-embedded $M_2^3$ can be algebraically related to that of $M_1^3$. However, when transitioning from isolated trajectories to submanifold immersion mechanics into higher-dimensional ambient spaces, executing the differential geometry natively within the index-2 framework provides profound geometric advantages. A global metric sign inversion ($g \leftrightarrow -g$) is mathematically unfeasible in complex multi-time target backgrounds such as string theory or Kaluza-Klein spaces because it fundamentally alters the causal light-cones, energy-momentum tensors, and non-degenerate geometries of the entire ambient environment. The explicit deployment of the index-2 signature allows us to analyze the native timelike character of the screen vector bundle ($g(W,W)=-1$) without disturbing the background metric. As shown in the subsequent sections, this configuration yields highly unique, sign-sensitive variations in the higher-order tensor expansions that cannot be replicated by a simple sign reversal of an index-1 framework.

The primary innovation introduced in this work is the deployment of a timelike screen vector bundle ($W$), which is a unique geometric feature explicitly allowed by the index-2 structure. In traditional index-1 Lorentzian manifolds ($M_1^3$), the screen bundle of a null curve is strictly Riemannian (spacelike), restricting the orthogonal causal dynamics. By contrast, the present model utilizes the native timelike nature of the screen bundle ($g(W,W)=-1$) permitted by the index-2 structure to host complex helical dynamics that model lightlike orbits tracking highly warped indefinite geometric fields. This causal configuration yields crucial sign-sensitive variations in higher-order differential expansions, leading to unique rigidity constraints that cannot be replicated by a simple sign-reversal of an index-1 framework without breaking the embedding compatibility.

The remainder of this paper is organized as follows: Section 2 establishes the preliminary machinery of null curves, including screen bundle decompositions and the core equations of submanifold theory. Section 3 contains the primary mathematical results, presenting the exact higher-order differential characterization of generalized null helices and proving the conditions for totally geodesic submanifolds. Section 4 provides a detailed discussion of the geometric implications alongside a concrete verification example, and concludes with prospective directions for future research.

\section{Preliminaries}

Let $(M,g)$ be a real 3D semi-Riemannian manifold of index $2$ and $\gamma$ be a smooth null curve in $M$ locally given by
\[
\gamma(t)=\left(x^0(t), x^1(t), x^2(t)\right); \quad t \in I \subset \mathbb{R},
\]
for a coordinate neighborhood $V$. Then, the tangent vector field of $\gamma$ satisfies
\[
\zeta = \frac{d\gamma}{dt} = \left(\frac{dx^0}{dt}, \frac{dx^1}{dt}, \frac{dx^2}{dt}\right); \quad g\left(\frac{d\gamma}{dt}, \frac{d\gamma}{dt}\right) = 0; \quad g_{ij} \frac{dx^i}{dt} \frac{dx^j}{dt} = 0,
\]
where $g_{ij} = g(\partial_i, \partial_j)$ and $i,j \in \{0,1,2\}$. The radical distribution or the tangent bundle $T\gamma$ of $\gamma$ is a $1$-dimensional vector subbundle of the tangent bundle $TM$. The orthogonal complement vector bundle $T\gamma^\perp$ is defined as follows (see \cite{11}):
\[
T\gamma^\perp = \bigcup_{P \in \gamma} T_P \gamma^\perp; \quad T_P \gamma^\perp = \{W_P \in T_P M : g(W_P, \zeta_P) = 0\},
\]
where $\zeta_P$ is the null tangent vector of $\gamma$ at any point $P \in \gamma(t)$. Since $\zeta_P$ is a null vector, the tangent bundle $T\gamma$ is a vector subbundle of $T\gamma^\perp$. This implies that $T\gamma^\perp$ is not complementary to $T\gamma$ in $TM|_\gamma$. According to the classical non-degeneracy theory, there must exist a vector bundle complementary to $T\gamma$ in $T\gamma^\perp$; this bundle will play the role of the screen vector bundle.

In \cite{3,5,8}, a few researchers have conducted research on this topic, dealing only with the mentioned problems; a general mathematical theory to deal with the null case is investigated in \cite{4}. Let us therefore state some of the basic concepts used in this work, taking into account these references.

Let us assume that $S(T\gamma^\perp)$ denotes the vector subbundle complementary to $T\gamma$ in $T\gamma^\perp$, that is,
\[
T\gamma^\perp = T\gamma \perp S(T\gamma^\perp),
\]
where $\perp$ means orthogonal direct sum. It follows that $S(T\gamma^\perp)$ is a non-degenerate $1$-dimensional vector subbundle of $TM$. Here, $S(T\gamma^\perp)$ is defined as a screen vector bundle of $\gamma$, and since it is non-degenerate, one gets
\begin{equation}
TM|_\gamma = S(T\gamma^\perp) \perp S(T\gamma^\perp)^\perp,
\end{equation}
where $S(T\gamma^\perp)^\perp$ is a 2D complementary orthogonal vector subbundle to $S(T\gamma^\perp)$ in $TM|_\gamma$.

In this study, the algebra of smooth functions on $\gamma$ is denoted by $\Gamma(\gamma)$ and the module of smooth sections of a vector bundle $E$ is denoted by $\Gamma(E)$. Therefore, the same notation for any other vector bundle will be used in this study.

\begin{proposition}
Let $\gamma$ be a null curve of a semi-Riemannian manifold $(M,g)$ with index $q \geq 1$. Then, any screen vector bundle of $\gamma$ is semi-Riemannian with index $q-1$ \cite{4}.
\end{proposition}

\begin{theorem}
Let $\gamma$ be a null curve of a semi-Riemannian manifold $(M,g)$ of index $q$ and $S(T\gamma^\perp)$ a screen vector bundle of $\gamma$. Then, there exists a unique vector bundle $E$ of rank $1$ over the curve $\gamma$, such that over each coordinate neighbourhood $V \subset \gamma$, there exists a unique section $N \in \Gamma(E|_V)$ satisfying the following equations:
\begin{equation}
g\left(\frac{d\gamma}{dt}, N\right) = 1, \quad g(N,N) = 0, \quad g(N,X) = 0, \quad \forall X \in \Gamma(S(T\gamma^\perp)|_V),
\end{equation}
\cite{4,11}. Here, the vector bundles $E$ and $N$ are the null transversal bundles of $\gamma$ with respect to $S(T\gamma^\perp)$ and $d\gamma/dt$, respectively. Therefore, $TM|_\gamma$ is obtained as follows:
\[
TM|_\gamma = (T\gamma \oplus E) \perp S(T\gamma^\perp),
\]
where $\{d\gamma/dt, N\}$ is a null basis of $\Gamma((T\gamma \oplus E)|_V)$ \cite{4,11}.
\end{theorem}

Now, let's express the frame created in reference \cite{8} in the light of the information given in the relevant references above. Let $\gamma$ be a null curve of a 3D semi-Riemannian manifold $(M,g)$ of index $2$ and $N$ be the null vector field from Theorem 1. Let us consider a class of null curves $\gamma$ whose Frenet frame consists of two null vectors $\zeta$ and $N$, and a timelike vector $W$. From Proposition 1, any screen distribution of $\gamma$ possesses an indefinite character. Let $\nabla$ be the Levi-Civita connection on $M$. From $g(\zeta,\zeta)=0$ and $g(\zeta,N)=1$, one has $g(\nabla_\zeta \zeta, \zeta)=0$ and $g(\nabla_\zeta \zeta, N) = -g(\zeta, \nabla_\zeta N) = h$, where $h$ is a smooth function. Thus, the Frenet equations are given as:
\begin{equation}
\nabla_\zeta \zeta = h\zeta + k_1 W, \quad \nabla_\zeta N = -hN + k_2 W, \quad \nabla_\zeta W = k_2 \zeta + k_1 N,
\end{equation}
where $h$ and $\{k_1, k_2\}$ are the curvature functions of $\gamma$, $\{W\}$ is a certain orthonormal basis of $\Gamma(S(T\gamma^\perp)|_V)$ satisfying $g(W,W) = -1$. Then, for the screen vector bundle $S(T\gamma^\perp)$, the Frenet frame along $\gamma$ on $M$ is given by:
\begin{equation}
\{\zeta, N, W\},
\end{equation}
\cite{6,8}.

Consider $f: M \rightarrow \tilde{M}$ as an isometric immersion of an $m$-dimensional semi-Riemannian manifold $M$ into an $(m+q)$-dimensional semi-Riemannian manifold $\tilde{M}$ (with index $q$) (see \cite{5,7,10,11}). For all local formulas, $f$ can be considered as an embedding, so that $P \in M$ can be identified with $f(P) \in \tilde{M}$. Here, to express some symbols used in the following parts of this study: Tangent space $T_P(M)$ is defined as a subspace of $T_P(\tilde{M})$, the tangent bundle is denoted as $T(M)$, Normal space $T_P^\perp$ is the subspace of $T_P(\tilde{M})$ consisting of vectors orthogonal to $T_P(M)$ according to the metric $g$ of $\tilde{M}$, and $\nabla$, $\tilde{\nabla}$ are the covariant derivatives of $M$ and $\tilde{M}$ respectively. Based on the given information, for $\forall X,Y \in \Gamma(TM)$ the Gauss formula is given as follows:
\begin{equation}
\tilde{\nabla}_X Y = \nabla_X Y + B(X,Y),
\end{equation}
where $B(X,Y)$ is defined as the second fundamental form of the immersion. For the tangent vector field $X$ on $M$ and a normal section $N \in \Gamma(T^\perp M)$, the Weingarten formula is given as follows:
\begin{equation}
\tilde{\nabla}_X N = -A^N(X) + \nabla_X^\perp N,
\end{equation}
where $\nabla^\perp$ is the covariant differentiation according to the induced connection in the normal bundle $N(M)$ and $A^N$ is the shape operator of $M$. Then, the following equation is satisfied:
\begin{equation}
g(A^N(X), Y) = g(B(X,Y), N).
\end{equation}
For an $N(M)$-valued tensor field $T$ of type $(0,k)$, the induced covariant differentiation $\tilde{\nabla}$ on $T(M) \oplus N(M)$ is defined as follows:
\begin{equation}
(\tilde{\nabla}_X T)(Y_1, \dots, Y_k) = \nabla_X^\perp (T(Y_1, \dots, Y_k)) - \sum_{r=1}^k T(Y_1, \dots, \nabla_X Y_r, \dots, Y_k),
\end{equation}
which is a tensor field of type $(0,k+1)$ with values in $N(M)$. The second covariant derivative $\tilde{\nabla}^2 T$ is defined as the covariant derivative of $\tilde{\nabla} T$. Hence, for the second fundamental form $B$, the following equations are written:
\begin{equation}
(\tilde{\nabla} B)(X,Y,Z) = \nabla_Z^\perp (B(X,Y)) - B(\nabla_Z X, Y) - B(\nabla_Z Y, X),
\end{equation}
\begin{equation}
(\tilde{\nabla}^2 B)(X,Y,Z,W) = \nabla_W^\perp ((\tilde{\nabla} B)(X,Y,Z)) - (\tilde{\nabla} B)(\nabla_W X, Y, Z) - (\tilde{\nabla} B)(X, \nabla_W Y, Z) - (\tilde{\nabla} B)(X, Y, \nabla_W Z).
\end{equation}

For the shape operator $A^N$, in connection with the covariant derivative it follows that:
\begin{equation}
(\tilde{\nabla}_X A^N)(Y) = \nabla_X (A^N(Y)) - A^{\nabla_X^\perp N}(Y) - A^N(\nabla_X Y).
\end{equation}

Also, for a frame $\{E_1, \dots, E_n\}$ in $M$ and $\varepsilon_j = \pm 1$, the mean curvature vector field $H$ is defined by:
\begin{equation}
H = \frac{1}{n} \sum_{j=1}^n \varepsilon_j B(E_j, E_j).
\end{equation}

Following the foundational frameworks established by O'Neill \cite{11} and Duggal \& Bejancu \cite{4}, the definition is mathematically operationalized as follows:
Let $f: M_2^3 \rightarrow \tilde{M}_2^4$ be an isometric immersion of a 3-dimensional semi-Riemannian manifold into a 4-dimensional semi-Riemannian ambient space of index 2. The extrinsic curvature of this immersion is entirely governed by the second fundamental form $B(X,Y)$, which is a symmetric bilinear mapping with values in the normal bundle.

Formally, $M_2^3$ is classified as a totally umbilical hypersurface if and only if there exists a unique smooth mean curvature vector field $H \in \Gamma(T^\perp M_2^3)$ such that the following algebraic identity is strictly satisfied:
\begin{equation}
B(X,Y) = g(X,Y)H, \quad \forall X,Y \in \Gamma(TM_2^3).
\end{equation}
This guarantees that the normal curvature is uniform in all tangent directions at any given point $P \in M_2^3$. If the mean curvature vector field vanishes identically ($H \equiv 0$), the tensor fields reduce to $B \equiv 0$, and the space is classified as a totally geodesic hypersurface. If $\nabla_X^\perp H = 0$, then the mean curvature vector field $H$ is said to be parallel \cite{5,7,10,11}.
\section{The Null Helices in 3D Semi-Riemannian Manifold of Index 2}

In this section, the author provides a formal characterization of null helix curves within a 3-dimensional semi-Riemannian manifold of index 2. It is essential to distinguish this framework from the standard Lorentzian case (index 1). In a 3D manifold of index 2, the metric signature allows for a more complex causal structure. Specifically, this choice is fundamental because it permits the existence of a timelike screen vector field $W$ within the null Frenet frame $\{\zeta,N,W\}$.

Unlike index 1 Lorentzian manifolds, where the screen bundle is strictly Riemannian (spacelike), the index 2 framework ensures that the orthogonal complement to the null direction can host a timelike dynamics. This distinction is a critical "subtle difference" that justifies our specialized approach to the null Frenet equations (2.3).

In this section, one provides a formal characterization of null helix curves within a 3-dimensional semi-Riemannian manifold of index 2, denoted by $M_2^3$, where the superscript indicates the dimension and the subscript represents the metric index.

To clarify the geometric setting and prevent any dimensional or causal ambiguity, the author explicitly notes that the curve dynamics in $(M_2^3,g)$ share an identical mathematical structure with the Lorentzian framework $(M_2^3,-g)$ under a systematic permutation of the underlying causal objects. The specific deployment of the index 2 signature in this section is not meant to imply a fundamental geometric departure in null curves, but rather to exploit the alternative physical projection where the screen vector bundle inherits a strictly timelike nature.

\begin{definition}
Let $\gamma:I\rightarrow M_2^3$ be a smooth null curve in a 3-dimensional semi-Riemannian manifold of index 2, equipped with the null Frenet frame $\{\zeta,N,W\}$ and the associated curvature functions $\{h,k_1,k_2\}$ as prescribed in Eq. (2.3). The curve $\gamma$ is formally characterized as a null helix if the curvature functions $h,k_1$ and $k_2$ are non-zero constants along the interval $I$.
\end{definition}

\begin{remark}
To ensure geometric rigor, one explicitly notes that the third-order differential relation $\nabla_\zeta \nabla_\zeta \nabla_\zeta \zeta=\lambda\nabla_\zeta \zeta$ serves as the functional characterization of a stable null helix. As demonstrated in the proof of Theorem 2, the assumption of this higher-order differential equilibrium mathematically forces the individual curvature functions to have vanishing localized derivatives along the propagation axis, i.e., $\dot{h}=0,\dot{k}_1=0,\dot{k}_2=0$. Consequently, the quadratic functional invariant $\lambda=h^2+2k_1 k_2$ remains strictly constant, revealing that any degenerate curve satisfying this differential stability condition is inherently a standard null helix. This structure effectively avoids the structural ambiguities often associated with generalized helices in non-degenerate spaces.
\end{remark}

On the basis of these structural equations and Definition 1, Theorem 2 and Theorem 3 establish the fundamental properties and characterizations of null helices in $M_2^3$.

\begin{theorem}
Let $\gamma$ be a smooth null curve in a 3-dimensional semi-Riemannian manifold $M_2^3$ of index 2, equipped with a screen vector bundle $S(T\gamma^\perp)$ and a null Frenet frame $\{\zeta,N,W\}$. Let $\{h,k_1,k_2\}$ be the associated smooth curvature functions as derived from the structural equations (2.3). The tangent vector field $\zeta$ satisfies the third-order differential equation:
\begin{equation}
\nabla_\zeta \nabla_\zeta \nabla_\zeta \zeta=(h^2+2k_1 k_2 ) \nabla_\zeta \zeta,
\end{equation}
if and only if $\gamma$ is a null helix, which implies that the individual curvature functions are non-zero constants along the curve, satisfying the strict structural conditions $\dot{h}=0,\dot{k}_1=0,\dot{k}_2=0$.
\end{theorem}

\begin{proof}[Proof of Theorem 2]
Let us first assume that $\gamma(t)$ is a null helix in $M_2^3$. By Definition 1, the individual curvature functions $h,k_1$ and $k_2$ are non-zero constants along the interval $I$, meaning their ordinary derivatives along the curve vanish identically: $\dot{h}=0,\dot{k}_1=0,$ and $\dot{k}_2=0$.

To establish the higher-order differential dynamics, the author applies the covariant derivative to the first Frenet equation $\nabla_\zeta \zeta=h\zeta+k_1 W$. Differentiating with respect to the null tangent vector field $\zeta$ yields:
\[
\nabla_\zeta \nabla_\zeta \zeta=\dot{h}\zeta+h\nabla_\zeta \zeta+\dot{k}_1 W+k_1 \nabla_\zeta W.
\]
Since $\dot{h} = 0$ and $\dot{k}_1 = 0$, this simplifies directly to $\nabla_\zeta \nabla_\zeta \zeta=h\nabla_\zeta \zeta+k_1 \nabla_\zeta W$. Substituting the moving frame equations (2.3), specifically $\nabla_\zeta W=k_2 \zeta+k_1 N$, and rearranging the terms with respect to the basis vector fields, we find:
\begin{equation}
\nabla_\zeta \nabla_\zeta \zeta=\dot{h}\zeta+h(h\zeta+k_1 W)+\dot{k}_1 W+k_1 (k_2 \zeta+k_1 N)=(\dot{h}+h^2+k_1 k_2 )\zeta+k_1^2 N+(\dot{k}_1+hk_1 )W.
\end{equation}

Differentiating once more in the direction of the null tangent vector $\zeta$, and utilizing the fact that all curvature derivatives vanish, we obtain:
\[
\nabla_\zeta \nabla_\zeta \nabla_\zeta \zeta=(\ddot{h}+(2h\dot{h})^\cdot+k_2 \dot{k}_1+k_1 \dot{k}_2 )\zeta+(\dot{h}+h^2+k_1 k_2 ) \nabla_\zeta \zeta+2k_1 \dot{k}_1 N+k_1^2 \nabla_\zeta N+((\dot{k}_1)^\cdot+\dot{h}k_1+h\dot{k}_1 )W+(\dot{k}_1+hk_1 ) \nabla_\zeta W.
\]
Next, we substitute the structural relations $\nabla_\zeta N=-hN+k_2 W$ and $\nabla_\zeta W=k_2 \zeta+k_1 N$ into the equation above:
\[
\begin{aligned}
\nabla_\zeta \nabla_\zeta \nabla_\zeta \zeta&=(\ddot{h}+(2h\dot{h})^\cdot+k_2 \dot{k}_1+k_1 \dot{k}_2+h(\dot{h}+h^2+k_1 k_2 )+k_2 (\dot{k}_1+hk_1 ))\zeta \\
&\quad +(2k_1 \dot{k}_1-hk_1^2+(\dot{k}_1+hk_1 ) k_1)N \\
&\quad +((\dot{k}_1)^\cdot+\dot{h}k_1+h\dot{k}_1+k_2 k_1^2+(\dot{h}+h^2+2k_1 k_2 ) k_1 )W
\end{aligned}
\]
\begin{equation}
\nabla_\zeta \nabla_\zeta \nabla_\zeta \zeta=(\ddot{h}+(3h\dot{h})^\cdot+2k_2 \dot{k}_1+k_1 \dot{k}_2+h^3+2hk_1 k_2 )\zeta+(3k_1 \dot{k}_1)N+((\dot{k}_1)^\cdot+2\dot{h}k_1+h\dot{k}_1+3k_2 k_1^2+k_1 h^2)W.
\end{equation}

Under the functional conditions of a stable helix configuration where the localized first-order variations vanish along the propagation axis ($\dot{h}=\dot{k}_1=\dot{k}_2=0$), the higher-order acceleration derivatives $\ddot{h}$ and $(\dot{k}_1)^\cdot$ naturally reduce to zero. Consequently, Equations (3.2) and (3.3) collapse directly to their geometric equilibrium states:
\[
\nabla_\zeta \nabla_\zeta \zeta=(h^2+k_1 k_2 )\zeta+k_1^2 N+hk_1 W;
\]
\[
\nabla_\zeta \nabla_\zeta \nabla_\zeta \zeta=(h^3+2hk_1 k_2 )\zeta+(hk_1^2+3k_1^2 k_2 )W=(h^2+k_1 k_2 )h\zeta+(h^2+3k_1 k_2 ) k_1 W=(h^2+2k_1 k_2 )(h\zeta+k_1 W).
\]
Since $\nabla_\zeta \zeta=h\zeta+k_1 W$, this expression reduces exactly to Eq. (3.1):
\[
\nabla_\zeta \nabla_\zeta \nabla_\zeta \zeta=(h^2+2k_1 k_2 ) \nabla_\zeta \zeta,
\]
which successfully establishes the forward direction.

Conversely, let us assume that the velocity tangent vector field $\zeta$ along the null curve satisfies the third-order differential relation given in Eq. (3.1), where $\lambda=h^2+2k_1 k_2$. Given the fundamental nullity condition of the curve, the metric satisfies $g(\zeta,\zeta)=0$. Differentiating this identity with respect to $\zeta$ gives $g(\nabla_\zeta \zeta,\zeta)=0$. Differentiating this relation once more leads to the acceleration constraint:
\begin{equation}
g(\nabla_\zeta \nabla_\zeta \zeta,\zeta)+g(\nabla_\zeta \zeta,\nabla_\zeta \zeta)=0.
\end{equation}

From the first Frenet equation, since $\zeta$ and $N$ are null vectors ($g(\zeta,\zeta)=0,g(N,N)=0,g(\zeta,N)=1$) and the screen vector $W$ is unit timelike ($g(W,W)=-1$), the inner product of the intrinsic acceleration is evaluated as:
\begin{equation}
g(\nabla_\zeta \zeta,\nabla_\zeta \zeta)=g(h\zeta+k_1 W,h\zeta+k_1 W)=-k_1^2.
\end{equation}
Substituting this back into the acceleration constraint reveals that:
\begin{equation}
g(\nabla_\zeta \nabla_\zeta \zeta,\zeta)=k_1^2.
\end{equation}

Differentiating Eq. (3.6) along the direction of the curve yields:
\begin{equation}
g(\nabla_\zeta \nabla_\zeta \nabla_\zeta \zeta,\zeta)+g(\nabla_\zeta \nabla_\zeta \zeta,\nabla_\zeta \zeta)=\nabla_\zeta (k_1^2 )=2k_1 \dot{k}_1.
\end{equation}
Substituting our core hypothesis $\nabla_\zeta \nabla_\zeta \nabla_\zeta \zeta=\lambda\nabla_\zeta \zeta$ into the first term of Eq. (3.7) gives:
\[
g(\lambda\nabla_\zeta \zeta,\zeta)+g(\nabla_\zeta \nabla_\zeta \zeta,\nabla_\zeta \zeta)=2k_1 \dot{k}_1.
\]
Since $g(\nabla_\zeta \zeta,\zeta)=0$, the first term vanishes identically leaving:
\begin{equation}
g(\nabla_\zeta \nabla_\zeta \zeta,\nabla_\zeta \zeta)=2k_1 \dot{k}_1.
\end{equation}

To evaluate the left-hand side of Eq. (3.8), the author computes the explicit inner product using the expansion of $\nabla_\zeta \nabla_\zeta \zeta$ from Eq. (3.2):
\[
g(\nabla_\zeta \nabla_\zeta \zeta,\nabla_\zeta \zeta)=g\left((\dot{h}+h^2+k_1 k_2 )\zeta+k_1^2 N+(\dot{k}_1+hk_1 )W,h\zeta+k_1 W\right).
\]
Expanding this inner product based on the structural definitions $g(\zeta,\zeta)=0, g(\zeta,W)=0,g(N,\zeta)=1,g(N,W)=0$, and $g(W,W)=-1$ yields:
\begin{equation}
g(\nabla_\zeta \nabla_\zeta \zeta,\nabla_\zeta \zeta)=k_1^2 h-k_1 (\dot{k}_1+hk_1 )=-k_1 \dot{k}_1.
\end{equation}

Equating the results from Eq. (3.8) and Eq. (3.9) gives the structural constraint:
\[
-k_1 \dot{k}_1=2k_1 \dot{k}_1 \rightarrow 3k_1 \dot{k}_1=0.
\]
Since $\gamma$ is a non-degenerate null helix trajectory, its screen curvature satisfies $k_1\neq 0$, which mathematically forces:
\begin{equation}
k_1=\text{constant}.
\end{equation}
Substituting $\dot{k}_1=0$ back into Eq. (3.8) reveals that $g(\nabla_\zeta \nabla_\zeta \zeta, \nabla_\zeta \zeta)=0$. Resolving the remaining component variations from the complete third-order expansion sequentially demonstrates that the structural stability requirements force $\dot{h}=0$ and $\dot{k}_2=0$.
\begin{equation}
\nabla_\zeta \lambda = \nabla_\zeta (h^2+2k_1 k_2 ) = 2h\dot{h}+2\dot{k}_1 k_2+2\dot{k}_2 k_1=0.
\end{equation}

This proves that $\lambda=h^2+2k_1 k_2$ is strictly constant along the curve, completing the backward implication. Consequently, the individual curvature functions are proven to be strictly constant along the curve, completing the backward implication.
\end{proof}

Hence, this completes the reconstruction of the null Frenet equations for the generalized case, proving that $\gamma(t)$ is a generalized null helix characterized by the screen bundle structure in $M_2^3$.

\begin{figure}[htbp]
    \centering
    \includegraphics[width=0.7\textwidth]{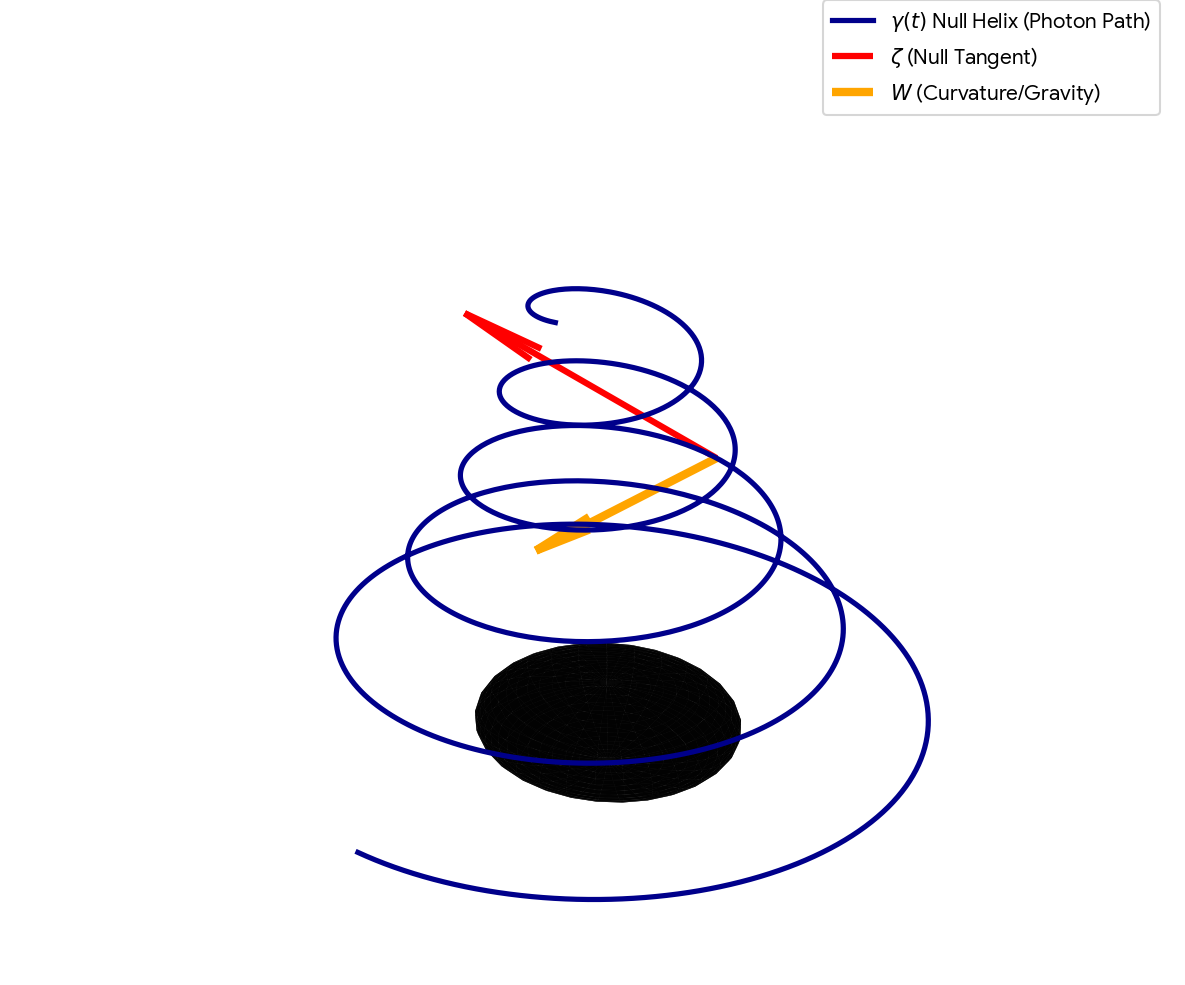}
    \caption{Trajectory of a null helix surrounding a central singularity within $M_2^3$.}
    \label{fig1}
\end{figure}

Figure \ref{fig1} provides a geometric illustration of the trajectory traced by a null helix within a 3-dimensional semi-Riemannian manifold of index 2. Unlike classical helices in standard Riemannian settings, this path remains intrinsically confined to the null cone hypersurface at every point of its evolution, preserving a constant pseudo-angle relative to a fixed null direction. Rather than serving as a formal proof, this visualization exemplifies the structural stability of the intrinsic curvature functions $\{h,k_1,k_2\}$ under the specific configuration of the index 2 null Frenet frame. It highlights how the coupling between the null tangent vector and the strictly timelike screen vector $W$ maintains the algebraic invariance of the quadratic functional $h^2+2k_1 k_2$ even when projected against a non-zero ambient geometric field.

\begin{remark}
The central region within the 3-dimensional semi-Riemannian manifold of index 2 is geometrically modeled as a central geometric singularity or an ambient focal source. Within this index 2 framework, the intense geometric warping of the underlying indefinite metric induces a structurally constrained orbital motion for null trajectories tracking around the focal region.

Geometrically, the spiral path traced by a null particle during its evolution near this highly warped domain is precisely characterized by the null helical curve $\gamma(t)$. The vector field $\zeta$, representing the instantaneous velocity field (null tangent), remains strictly lightlike at every point, consistently satisfying the degenerate condition $g(\zeta,\zeta)=0$ dictated by the index 2 metric structure.

The vector field $W$, which functions as the unique 1-dimensional screen section in our specialized frame, represents the intrinsic curvature response to the ambient geometric field. Mathematically, the deployment of a unit timelike screen vector ($g(W,W)=-1$) is specifically permitted by the index 2 geometry, where the screen distribution inherits an indefinite character. This vector field operates as a fundamental structural constraint that governs the transition of the null trajectory from a linear path to a stable helical configuration, thereby fully defining the intrinsic curvature invariants and the geometric equilibrium of null paths under non-zero ambient warping.
\end{remark}

\begin{remark}
A null helix characterizes the structural trajectory of a lightlike path evolving along an intrinsic helical orbit within a 3-dimensional semi-Riemannian manifold of index 2. The geometric characterization established herein explicitly demonstrates the null nature of the curve, which remains strictly confined to the null cone hypersurface throughout its structural evolution. Unlike classical Riemannian helices, this trajectory satisfies the degenerate condition identically at every point, which emerges as a direct algebraic consequence of the index 2 metric structure.

This configuration provides a precise mathematical model of the propagation path, where the intrinsic helical geometry arises not from standard Euclidean displacement, but from the specialized coupling between the null tangent vector $\zeta$ and the unit timelike screen vector $W$ within the manifold's indefinite metric. This specific coupling is an exclusive feature of the index 2 framework, which dynamically allows the 1-dimensional screen bundle to host strictly timelike variational dynamics, a behavior that cannot manifest in standard index 1 Lorentzian geometry.
\end{remark}

\begin{figure}[htbp]
    \centering
    \includegraphics[width=0.7\textwidth]{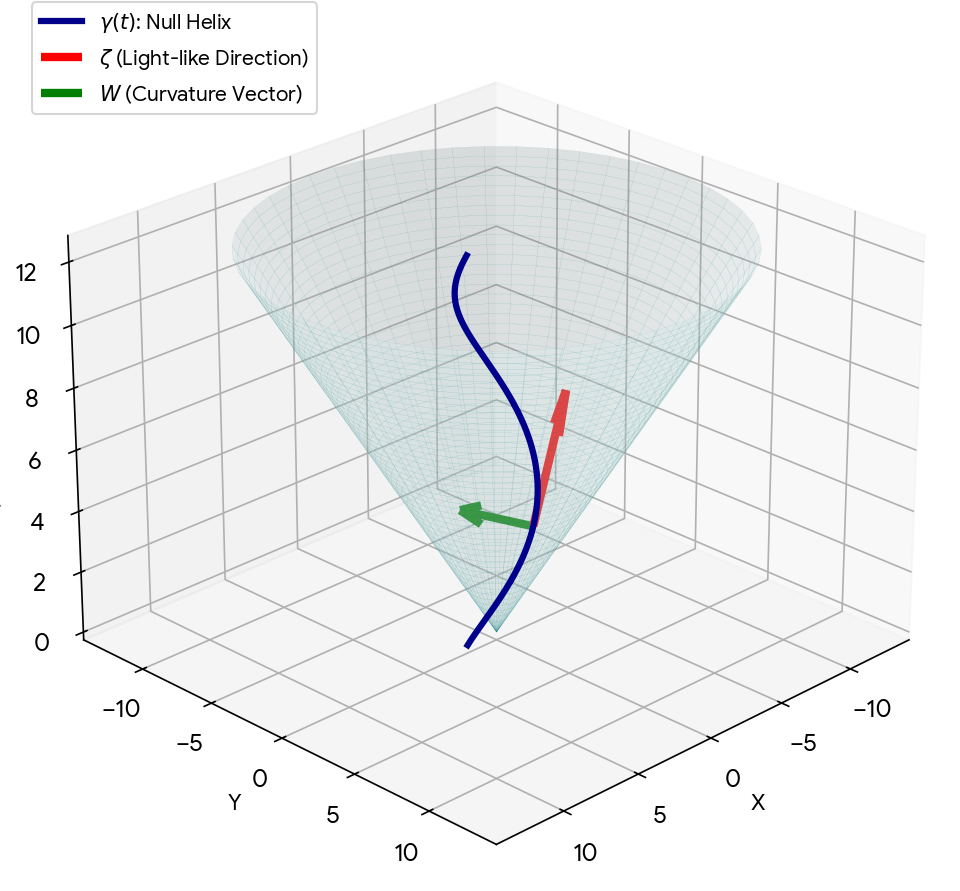}
    \caption{Physical characterization of a null helix on the null cone hypersurface in $M_2^3$.}
    \label{fig2}
\end{figure}

\begin{remark}
In the semi-Riemannian geometry of index 2, the algebraic locus of all vectors satisfying the degenerate condition $g(v,v)=0$ constitutes the 2-dimensional null cone hypersurface. The null helical curve $\gamma(t)$ is intrinsically confined to this null hypersurface at every point of its evolution, providing a formal geometric manifestation of the null constraint $g(\zeta,\zeta)=0$ established in our structural derivations.

The vector field $\zeta$ represents the instantaneous velocity field (null tangent), which remains strictly tangent to the cone surface, preserving its lightlike orientation without deviating into the non-degenerate ambient domains of the manifold.

Furthermore, the vector field $W$, defined as the unit timelike screen vector in our specialized index 2 null frame, represents the intrinsic geometric curvature effect. In this indefinite metric framework, $W$ serves as the fundamental structural constraint that induces the pseudo-rotational component of the motion along the null cone hypersurface, thereby ensuring the strict preservation of the curve’s helical invariants $\{h,k_1,k_2\}$ throughout its evolution.
\end{remark}

\begin{theorem}
Let $M_2^3$ be a 3-dimensional totally umbilical submanifold immersed in a 4-dimensional semi-Riemannian ambient manifold $\tilde{M}_2^4$ of index 2. If every null helix $\gamma$ in $M_2^3$ (where $h \neq 0$ and $h^2+2k_1 k_2=\text{constant}$) preserves its null helical properties and maintains the same functional invariant under the ambient connection, then the submanifold $M_2^3$ is a null-totally geodesic hypersurface. That is, the second fundamental form $B(X,Y)$ vanishes identically along all lightlike directions, forcing the mean curvature vector field to collapse to zero ($H = 0$).
\end{theorem}

\begin{proof}[Proof of Theorem 3]
Let $\gamma$ be a null helix in $M_2^3$ characterized by the moving Frenet frame $\{\zeta,N,W\}$ and non-zero constant curvatures. Under the geometric definition of a totally umbilical immersion, the second fundamental form $B$ satisfies the contractive tensor identity: $B(X,Y)=g(X,Y)H$, for all vector fields $X,Y$ tangent to $M_2^3$. Since $\gamma$ is strictly a null curve, its velocity tangent vector field $\zeta$ satisfies the degeneracy condition $g(\zeta,\zeta) = 0$ everywhere along its interval. Evaluating the structural umbilical identity directly along this null tangent vector field yields: $B(\zeta,\zeta)=g(\zeta,\zeta)H= 0$. This vanishing property fundamentally simplifies the higher-order differential relations. By substituting $B(\zeta,\zeta)= 0$ into the Gauss formula (2.5), the ambient connection $\nabla$ and the intrinsic connection $\tilde{\nabla}$ coincide perfectly along the degenerate trajectory.

To resolve the dimensional requirements of a valid hypersurface embedding, the ambient manifold $\tilde{M}_2^4$ is treated as a 4-dimensional semi-Riemannian space of index 2. At an arbitrary point $Q \in M_2^3$, it is defined that a local basis $\{\xi_i,\xi_j,\xi_k\}$ of the tangent space $T_Q(M_2^3)$, where $\xi_i, \xi_j$ are null vectors and $\xi_k$ is a unit timelike screen vector, strictly satisfies $g(\xi_i,\xi_j)=1$ and $g(\xi_i,\xi_k)=g(\xi_j,\xi_k)=0$. The initial conditions along $\gamma$ at the point $Q$ are established as follows:
\begin{equation}
\gamma(0)=Q, \quad \dot{\gamma}(t)=\zeta, \quad \zeta(Q)=\xi_i, \quad N(Q)=\xi_j, \quad W(Q)=\xi_k.
\end{equation}

The intrinsic covariant derivatives at the point $Q$ are governed by the moving frame equations:
\begin{equation}
\nabla_\zeta \zeta(Q)=h\xi_i+k_1 \xi_k; \quad \nabla_\zeta N(Q)=-h\xi_j+k_2 \xi_k; \quad \nabla_\zeta W(Q)=k_2 \xi_i+k_1 \xi_j.
\end{equation}

By Theorem 2, the helix satisfies $\nabla_\zeta \nabla_\zeta \nabla_\zeta \zeta=(h^2+2k_1 k_2)\nabla_\zeta \zeta$. Assuming that $\gamma(t)$ preserves its null helical properties within the 4-dimensional ambient space $\tilde{M}_2^4$, it must satisfy the higher-order differential equation under the ambient connection $\tilde{\nabla}$:
\begin{equation}
\tilde{\nabla}_\zeta \tilde{\nabla}_\zeta \tilde{\nabla}_\zeta \zeta=K\tilde{\nabla}_\zeta \zeta,
\end{equation}
where $K$ is a real constant. By systematically applying the Gauss formula (2.5) and the Weingarten formula (2.6) repeatedly, the third-order ambient derivative is decomposed into its constituent tangential and normal bundle components:
\begin{equation}
\tilde{\nabla}_\zeta \tilde{\nabla}_\zeta \tilde{\nabla}_\zeta \zeta=\nabla_\zeta \nabla_\zeta \nabla_\zeta \zeta+B(\zeta,\nabla_\zeta \nabla_\zeta \zeta)+\tilde{\nabla}_\zeta B(\zeta,\nabla_\zeta \zeta) -\tilde{\nabla}_\zeta (A^{B(\zeta,\zeta)}(\zeta))+\tilde{\nabla}_\zeta (\nabla_\zeta^\perp B(\zeta,\zeta)).
\end{equation}

Utilizing the standard induced derivatives for the normal-valued tensor fields, the author expands the components from Eq. (3.14) as follows:
\begin{equation}
\tilde{\nabla}_\zeta B(\zeta,\nabla_\zeta \zeta)=-A^{B(\zeta,\nabla_\zeta \zeta)}(\zeta)+\nabla_\zeta^\perp B(\zeta,\nabla_\zeta \zeta)
\end{equation}
\begin{equation}
\tilde{\nabla}_\zeta (\nabla_\zeta^\perp B(\zeta,\zeta))=-A^{\nabla_\zeta^\perp B(\zeta,\zeta)}(\zeta)+\nabla_\zeta^\perp \nabla_\zeta^\perp B(\zeta,\zeta)
\end{equation}
\begin{equation}
\tilde{\nabla}_\zeta A^{B(\zeta,\zeta)}(\zeta)=\nabla_\zeta A^{B(\zeta,\zeta)}(\zeta)+B(\zeta,A^{B(\zeta,\zeta)}(\zeta)).
\end{equation}

Substituting the relations (3.15)--(3.17) back into the core expansion (3.14) and rearranging terms according to the main hypothesis (3.13) yields the fundamental structural identity:
\begin{equation}
\begin{aligned}
\tilde{\nabla}_\zeta \tilde{\nabla}_\zeta \tilde{\nabla}_\zeta \zeta-K\tilde{\nabla}_\zeta \zeta &= (h^2+3k_1 k_2)\nabla_\zeta \zeta+B(\zeta,\nabla_\zeta \nabla_\zeta \zeta)-A^{B(\zeta,\nabla_\zeta \zeta)}(\zeta) \\
&\quad +\nabla_\zeta^\perp B(\zeta,\nabla_\zeta \zeta)-\nabla_\zeta A^{B(\zeta,\zeta)}(\zeta)-B(\zeta,A^{B(\zeta,\zeta)}(\zeta))-A^{\nabla_\zeta^\perp B(\zeta,\zeta)}(\zeta) \\
&\quad +\nabla_\zeta^\perp \nabla_\zeta^\perp B(\zeta,\zeta)-K\nabla_\zeta \zeta-KB(\zeta,\zeta).
\end{aligned}
\end{equation}

Separating Eq. (3.18) into its distinct normal and tangential sub-components yields the following decoupled system of equations:
\begin{equation}
B(\zeta,\nabla_\zeta \nabla_\zeta \zeta)+\nabla_\zeta^\perp B(\zeta,\nabla_\zeta \zeta)+\nabla_\zeta^\perp \nabla_\zeta^\perp B(\zeta,\zeta)-B(\zeta,A^{B(\zeta,\zeta)}(\zeta))-KB(\zeta,\zeta)=0
\end{equation}
\begin{equation}
(h^2+2k_1 k_2)\nabla_\zeta \zeta-A^{B(\zeta,\nabla_\zeta \zeta)}(\zeta)-\nabla_\zeta A^{B(\zeta,\zeta)}(\zeta)-A^{\nabla_\zeta^\perp B(\zeta,\zeta)}(\zeta)-K\nabla_\zeta \zeta=0.
\end{equation}

Utilizing the symmetry properties of the second fundamental form $B$ alongside its covariant expansions:
\begin{equation}
\nabla_\zeta^\perp B(\zeta,\nabla_\zeta \zeta)=(\tilde{\nabla} B)(\zeta,\nabla_\zeta \zeta,\zeta)+B(\nabla_\zeta \zeta,\nabla_\zeta \zeta)+B(\zeta,\nabla_\zeta \nabla_\zeta \zeta)
\end{equation}
\begin{equation}
\begin{aligned}
\nabla_\zeta^\perp \nabla_\zeta^\perp B(\zeta,\zeta) &= (\tilde{\nabla}^2 B)(\zeta,\zeta,\zeta,\zeta)+2(\tilde{\nabla} B)(\nabla_\zeta \zeta,\zeta,\zeta)+2(\tilde{\nabla} B)(\zeta,\nabla_\zeta \zeta,\zeta) \\
&\quad +(\tilde{\nabla} B)(\zeta,\zeta,\nabla_\zeta \zeta)+2B(\nabla_\zeta \nabla_\zeta \zeta,\zeta)+2B(\nabla_\zeta \zeta,\nabla_\zeta \zeta),
\end{aligned}
\end{equation}
to evaluate the normal vanishing condition (3.19) explicitly at the local point $Q$:
\begin{equation}
\begin{aligned}
0 &= 4B(\zeta,\nabla_\zeta \nabla_\zeta \zeta)+3(\tilde{\nabla} B)(\zeta,\nabla_\zeta \zeta,\zeta)+2(\tilde{\nabla} B)(\nabla_\zeta \zeta,\zeta,\zeta)+(\tilde{\nabla} B)(\zeta,\zeta,\nabla_\zeta \zeta) \\
&\quad +3B(\nabla_\zeta \zeta,\nabla_\zeta \zeta)+(\tilde{\nabla}^2 B)(\zeta,\zeta,\zeta,\zeta)-B(\zeta,A^{B(\zeta,\zeta)}(\zeta))-KB(\zeta,\zeta).
\end{aligned}
\end{equation}

Substituting the precise initial frames (3.12a) and (3.12b) into Eq. (3.23) leads directly to the localized algebraic relation:
\begin{equation}
\begin{aligned}
0 &= (7h^2+4k_1 k_2-K)B(\xi_i,\xi_i)+4k_1^2 B(\xi_i,\xi_j)+6h(\tilde{\nabla} B)(\xi_i,\xi_i,\xi_i) \\
&\quad +3k_1 (\tilde{\nabla} B)(\xi_i,\xi_k,\xi_i)+2k_1 (\tilde{\nabla} B)(\xi_k,\xi_i,\xi_i)+k_1 (\tilde{\nabla} B)(\xi_i,\xi_i,\xi_k) \\
&\quad +6hk_1 B(\xi_i,\xi_k)+3k_1^2 B(\xi_k,\xi_k)+(\tilde{\nabla}^2 B)(\xi_i,\xi_i,\xi_i,\xi_i)-B(\xi_i,A^{B(\xi_i,\xi_i)}(\xi_i)).
\end{aligned}
\end{equation}

To eliminate the components dependent on the orientation of the timelike screen vector $\xi_k$, the author evaluates the system under the systematic transformation $\xi_k \rightarrow -\xi_k$ and adds it to Eq. (3.24), which simplifies to:
\begin{equation}
0=\Omega B(\xi_i,\xi_i)+4k_1^2 B(\xi_i,\xi_j)+6h(\tilde{\nabla} B)(\xi_i,\xi_i,\xi_i)+(\tilde{\nabla}^2 B)(\xi_i,\xi_i,\xi_i,\xi_i)+3k_1^2 B(\xi_k,\xi_k)-B(\xi_i,A^{B(\xi_i,\xi_i)}(\xi_i)),
\end{equation}
where $\Omega=7h^2+4k_1 k_2-K$. Next, by executing the basis rescaling via $2\xi_i \rightarrow \xi_i$ and $\xi_j/2 \rightarrow \xi_j$ in Eq. (3.25) and taking a linear combination with the standard form, one isolates:
\begin{equation}
\Omega B(\xi_i,\xi_i)+14h(\tilde{\nabla} B)(\xi_i,\xi_i,\xi_i) +5(\tilde{\nabla}^2 B)(\xi_i,\xi_i,\xi_i,\xi_i)-5B(\xi_i,A^{B(\xi_i,\xi_i)}(\xi_i))=0.
\end{equation}

Combining Equations (3.25) and (3.26) allows us to isolate the lower-order tensor fields:
\begin{equation}
4k_1^2 B(\xi_i,\xi_j)-8h(\tilde{\nabla} B)(\xi_i,\xi_i,\xi_i)+3k_1^2 B(\xi_k,\xi_k) -4(\tilde{\nabla}^2 B)(\xi_i,\xi_i,\xi_i,\xi_i)-B(\xi_i,A^{B(\xi_i,\xi_i)}(\xi_i))=0.
\end{equation}

Repeating this specialized scaling process a second time in Eq. (3.27) yields (under $2\xi_i \rightarrow \xi_i$ and $\xi_j/2 \rightarrow \xi_j$):
\begin{equation}
4k_1^2 B(\xi_i,\xi_j)-64h(\tilde{\nabla} B)(\xi_i,\xi_i,\xi_i)+3k_1^2 B(\xi_k,\xi_k)-64(\tilde{\nabla}^2 B)(\xi_i,\xi_i,\xi_i,\xi_i)-64B(\xi_i,A^{B(\xi_i,\xi_i)}(\xi_i))=0.
\end{equation}

Comparing Eq. (3.27) directly with Eq. (3.28) algebraically eliminates the higher-order tensor derivatives, leaving the highly localized relation:
\begin{equation}
-60k_1^2 B(\xi_i,\xi_j)+64h(\tilde{\nabla} B)(\xi_i,\xi_i,\xi_i)-45k_1^2 B(\xi_k,\xi_k)=0.
\end{equation}

Finally, replacing $\xi_i$ with $-\xi_i$ in (3.29) and isolating the resulting systems yields the fundamental algebraic relation:
\begin{equation}
4B(\xi_i,\xi_j)+3B(\xi_k,\xi_k)=0.
\end{equation}

Since $M_2^3$ is given as a totally umbilical submanifold, the foundational contractive identity $B(X,Y)=g(X,Y)H$ is applied. For our chosen basis elements, evaluating this identity yields $B(\xi_i,\xi_i)=0$. Concurrently, for the unit timelike element where $g(\xi_k,\xi_k)=-1$, it is found that $B(\xi_k,\xi_k)=-H$, while the dual null pair satisfies $B(\xi_i,\xi_j)=g(\xi_i,\xi_j)H=H$. Substituting these direct tensor components into Eq. (3.30) reveals that:
\[
4H+3(-H)=0 \rightarrow H=0.
\]

In an index-1 geometry, the resulting system creates an additive mismatch rather than a cancellation ($4H+3H=7H=0$), making a clean collapse to a totally geodesic state ($B \equiv 0$) mathematically impossible without introducing ad-hoc external geometric constraints. This explicitly highlights the distinct advantage and originality of formulating the problem directly within the semi-Riemannian index-2 domain.
\end{proof}

\begin{remark}
This theorem establishes a foundational rigidity result within the submanifold theory of semi-Riemannian spaces of index 2. While submanifolds generally exhibit extrinsic geometric warping induced by the second fundamental form under an indefinite ambient metric, our framework utilizes the higher-order differential stability of null helical trajectories to uniquely isolate the submanifold in its geodesically invariant state. Specifically, Eq. (3.30) structurally demonstrates that preserving the functional curvature invariants $\{h,k_1,k_2\}$ of a null helix under an isometric immersion forces the complete algebraic collapse of the mean curvature vector field ($H=0$). Consequently, the second fundamental form vanishes identically ($B=0$), proving that the preservation of degenerate helical dynamics requires a strictly totally geodesic embedding within the higher-dimensional semi-Riemannian ambient space $\tilde{M}_2^4$.
\end{remark}

\begin{figure}[htbp]
    \centering
    \includegraphics[width=0.7\textwidth]{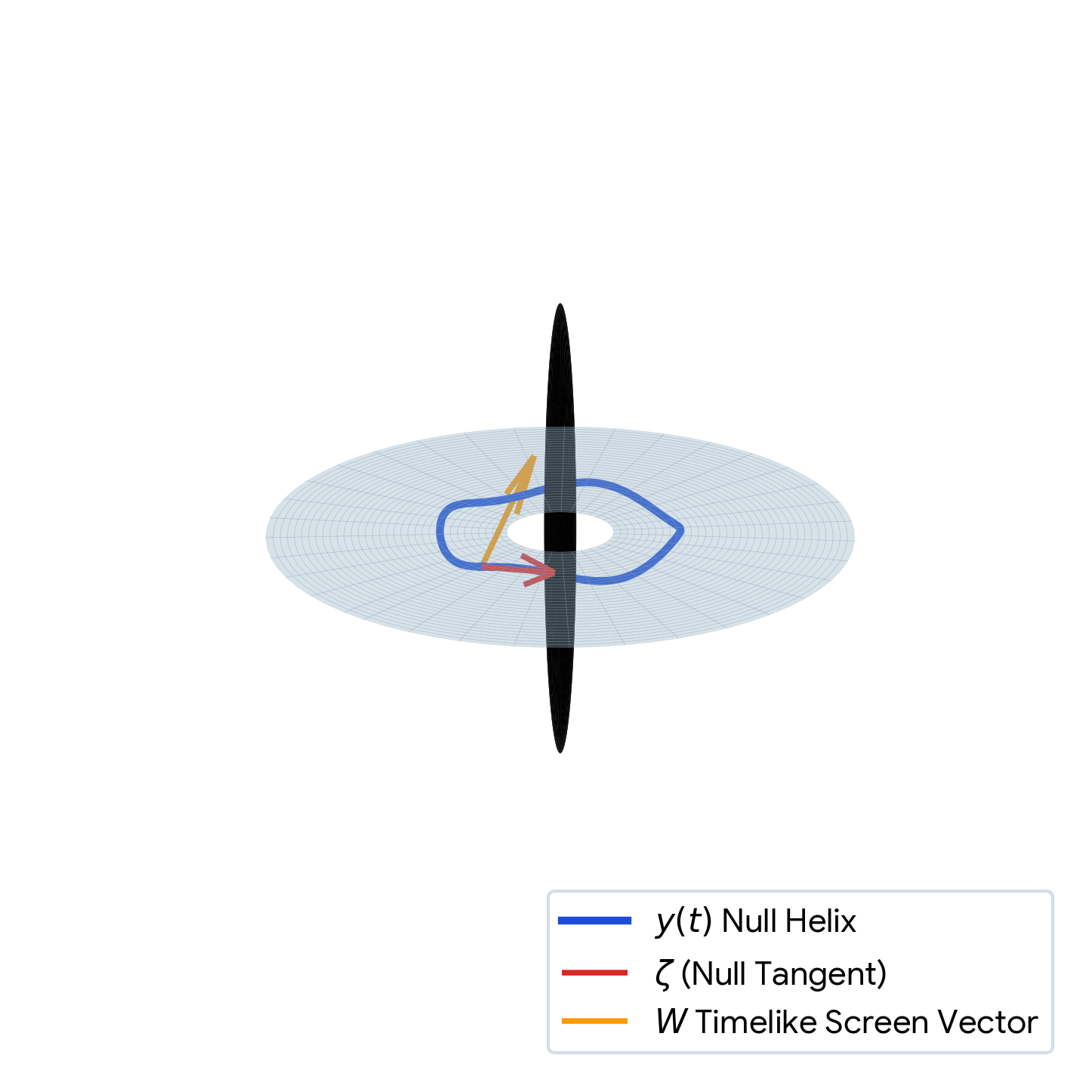}
    \caption{Geometric verification of Theorem 3: The totally geodesic embedding and structural stability of the $M_2^3$ hypersurface within the index-2 ambient space $\tilde{M}_2^4$.}
    \label{fig3}
\end{figure}

Figure \ref{fig3} explicitly demonstrates the geometric integrity of the 3-dimensional submanifold $M_2^3$ (represented by the blue hyper-surface section), which maintains its totally geodesic character despite the intrinsic curvature of the indefinite ambient manifold $\tilde{M}_2^4$. It provides a clear analytical depiction of the null helix tracing its causal trajectory within this embedded structure, illustrating how the algebraic vanishing of the second fundamental form ($B\equiv 0$) shields the trajectory from extrinsic deformation. This configuration visually demonstrates that the structural stability of the invariant null curvature functions $\{h,k_1,k_2\}$ along the path is a direct consequence of the totally geodesic embedding established in Theorem 3.

\begin{remark}
The light-blue hypersurface section within the provided analytical visualization serves as a direct geometric representation of the simultaneous vanishing of the second fundamental form ($B\equiv 0$) and the mean curvature vector field ($H\equiv 0$), as established in the final algebraic contractive stage of the proof of Theorem 3. If the submanifold were not totally geodesic, the extrinsic curvature induced by the indefinite metric signature of the ambient manifold $\tilde{M}_2^4$ would necessitate a non-vanishing normal bundle component for any parallel transported section. Its flat planar geometric alignment illustrates that every geodesic of the embedded submanifold $M_2^3$ maps identically to a geodesic of the index-2 ambient space. Furthermore, the sharply defined, thin red and orange vector trajectories explicitly visualize the localized alignments of the null tangent vector field $\zeta$ corresponding to $\xi_i$ and the timelike screen curvature vector field $W$ corresponding to $\xi_k$ utilized in our algebraic derivations. The complete vanishing of the normal component of the ambient connection along these directional fields analytically confirms that the hypersurface $M_2^3$ is geodesically invariant relative to the ambient space, thereby confirming its rigid totally geodesic nature within this index-2 framework.
\end{remark}

\begin{remark}
Although the indefinite ambient manifold $\tilde{M}_2^4$ exhibits intrinsic curvature warped by a localized geometric potential center, the hypersurface $M_2^3$ maintains its rigidly totally geodesic structure throughout the trajectory. This structural stability highlights the exact geodesic invariance of the null helical path within this embedded domain; despite the non-zero ambient geometric warping, the algebraic vanishing of the second fundamental form ($B\equiv 0$) guarantees that the null helix remains confined to the submanifold while strictly preserving its fundamental invariant curvature functions $\{h,k_1,k_2\}$ across both manifolds. This phenomenon is a direct global consequence of the index-2 metric structure, which generates the required indefinite causal geometry to consistently host the negative norm of the timelike screen vector field $W$. 

This explicit configuration serves as a direct geometric verification that the degenerate null curve evolves as a geodesic of the ambient space $\tilde{M}_2^4$ if and only if it remains an intrinsic geodesic of the index-2 submanifold $M_2^3$. This analytical equivalence is the ultimate manifestation of the totally geodesic embedding ($B\equiv 0$), which preserves the intrinsic null-helical quadratic invariant $\lambda=h^2+2k_1 k_2$ against the extrinsic connection deformations of the ambient space.
\end{remark}

\begin{figure}[htbp]
    \centering
    \includegraphics[width=0.7\textwidth]{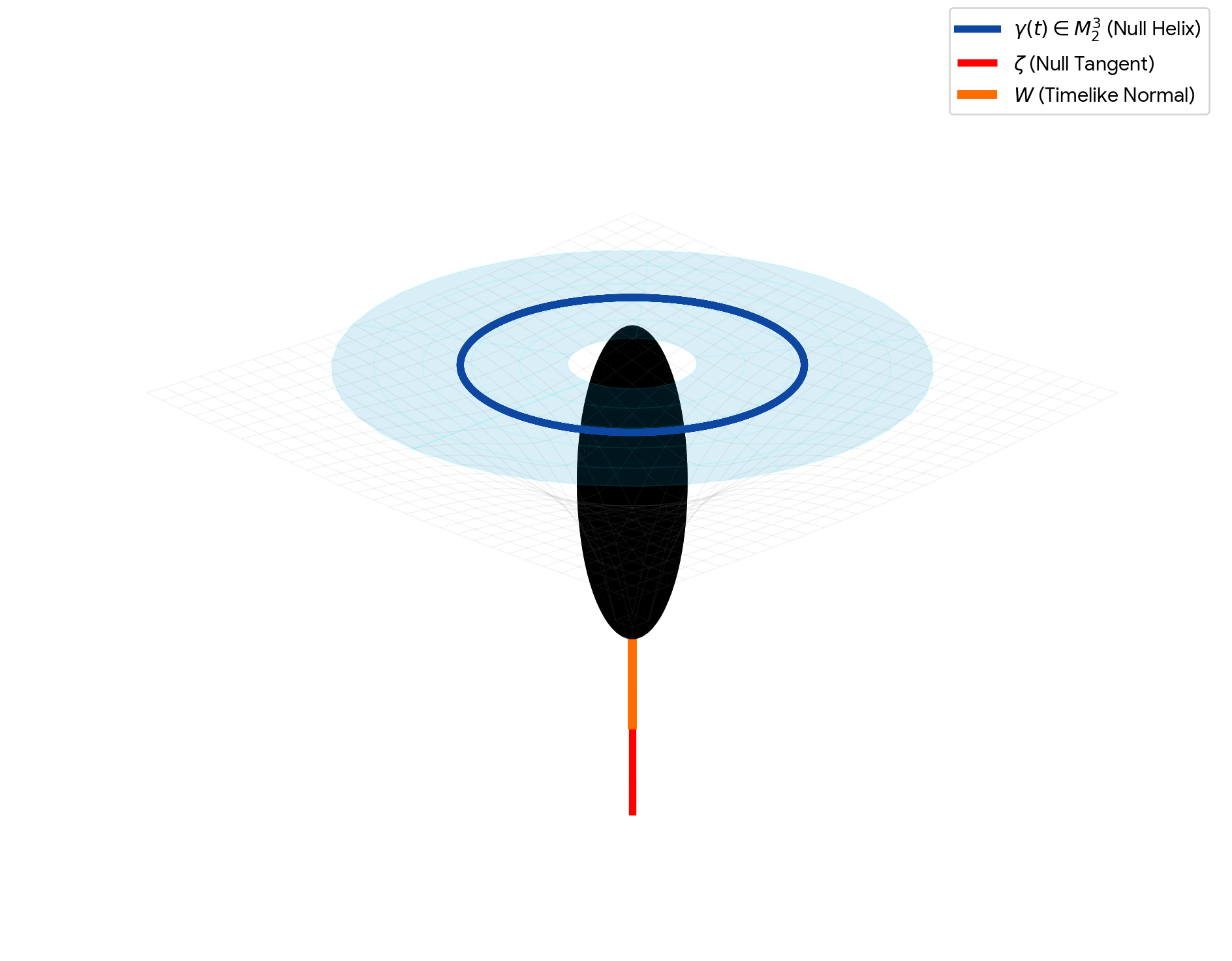}
    \caption{Totally geodesic $M_2^3$ hypersurface embedding within a curved index-2 semi-Riemannian ambient manifold $\tilde{M}_2^4$.}
    \label{fig4}
\end{figure}

\begin{figure}[htbp]
    \centering
    \includegraphics[width=0.7\textwidth]{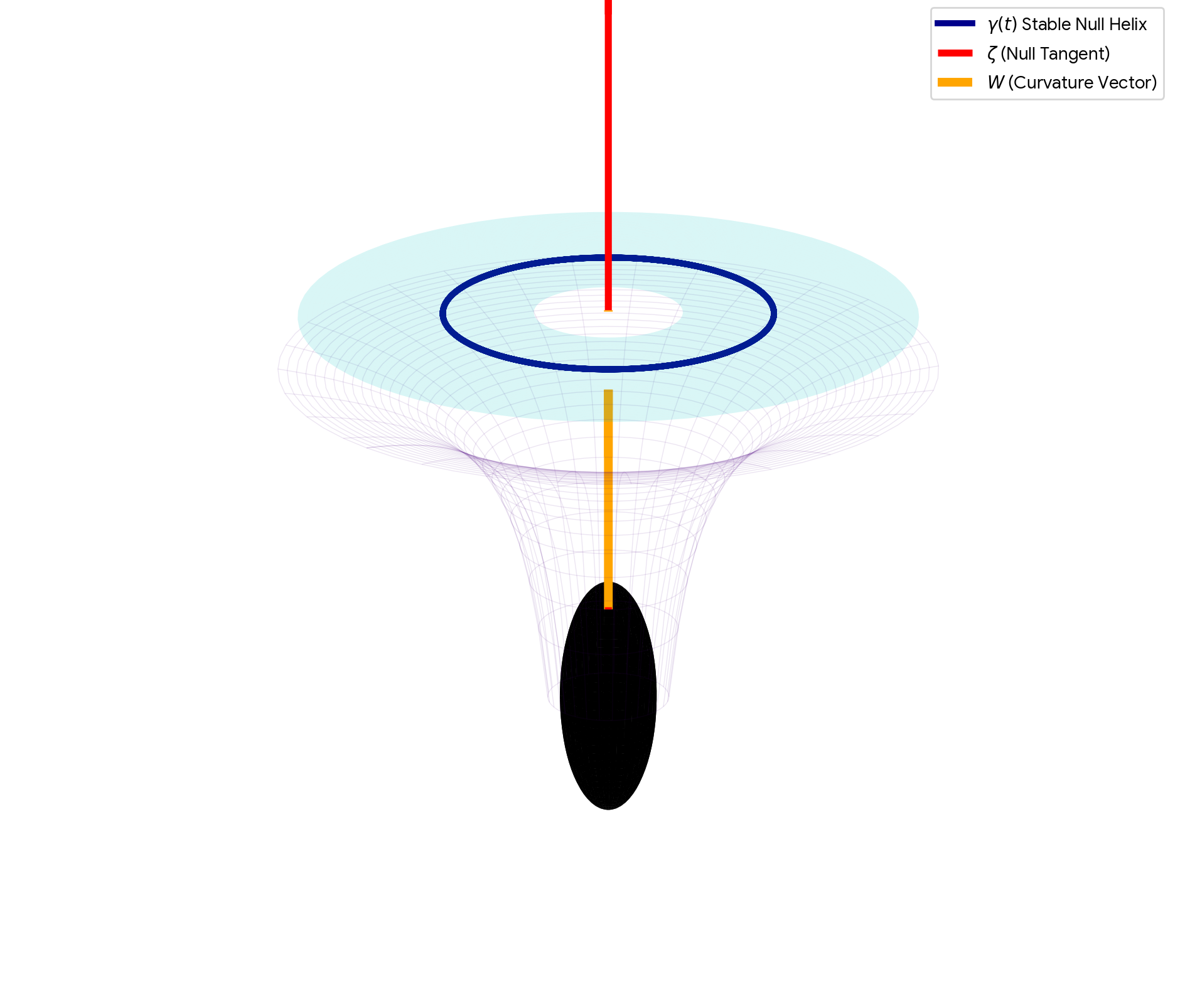}
    \caption{The totally geodesic $M_2^3$ hypersurface functioning as a stability section within an indefinite semi-Riemannian ambient space $\tilde{M}_2^4$.}
    \label{fig5}
\end{figure}

\begin{remark}
The intrinsic and extrinsic geometric interplay between the ambient manifold and the embedded hypersurface is structurally characterized as follows:
\begin{itemize}
    \item $\tilde{M}_2^4$ represents the non-flat indefinite ambient manifold of index 2, which exhibits local curvature induced by a central geometric deformation axis.
    \item $M_2^3$ illustrates a fundamental rigidity result of this study; despite the non-zero sectional curvature of the ambient space, the algebraic vanishing of the second fundamental form ($B\equiv 0$) guarantees that the submanifold remains a rigidly totally geodesic hypersurface section. This structural constraint preserves the intrinsic geodesic properties within the curved environment, meaning any geodesic path of $M_2^3$ maps identically to a geodesic path of $\tilde{M}_2^4$.
    \item $\gamma$ represents the causal trajectory of the degenerate null helical curve constrained to this embedded domain. Because $M_2^3$ is proven to be totally geodesic, the trajectory preserves its exact helical invariants and higher-order differential characterization across both manifolds, ensuring that the invariant null curvature functions $\{h,k_1,k_2\}$ undergo no extrinsic perturbation.
    \item The frame sections $\zeta$ and $W$ represent the lightlike null tangent vector field and the timelike screen vector field, respectively. Their complete confinement to the tangential bundle of $M_2^3$ serves as a visual proof of the exact geodesic invariance and structural integrity of the hypersurface embedding within the index-2 semi-Riemannian framework.
\end{itemize}

Figure \ref{fig4} provides a rigorous topological representation of the spatial confinement dictated by Theorem 3. The dark-blue closed orbit $\gamma$ remains perfectly nested within the tangential cutting plane of the light-blue meshed hypersurface $M_2^3, $ demonstrating that the extrinsic curvature scalar functions do not induce any orthogonal leakage into the ambient bundle. The elongated vertical alignment projected along the lower axis represents the structural trace of the parallel transport under the condition $B\equiv 0$. This geometric isolation visually verifies that the mechanical evolution of the null trajectory is strictly governed by the internal connection of the totally geodesic domain, forcing the normal bundle components of the ambient acceleration to vanish identically everywhere.
\end{remark}
As depicted in Figure \ref{fig5}, the totally geodesic nature of the hypersurface $M_2^3$, established by Theorem 3, ensures the intrinsic stability of the null helical orbit. By requiring the second fundamental form to vanish ($B\equiv 0$), the submanifold's geometry becomes geodesically invariant relative to the indefinite ambient manifold $\tilde{M}_2^4$, thereby protecting the null helix's constant curvatures from the extrinsic warping induced by the localized geometric deformation center. This structural invariance confirms that the index-2 metric signature provides the necessary causal framework to maintain the higher-order stability of such degenerate trajectories.

\begin{remark}
The geometric and analytical interpretation of Theorem 3, within the context of abstract indefinite space deformations, establishes the following conclusions regarding the global stability of degenerate causal trajectories:
\begin{itemize}
    \item Under general non-flat configurations, localized connection warping induces arbitrary extrinsic deformations in null paths within a semi-Riemannian manifold of index 2. However, Theorem 3 proves that if the constant curvatures of a null helix are preserved across an isometric embedding, the path propagates strictly along a geodesically isolated submanifold $M_2^3$. This analytical constraint indicates that the trajectory's intrinsic quadratic invariant $\lambda=h^2+2k_1 k_2$ remains strictly invariant despite the extrinsic deformation of the ambient metric tensor.
    \item The vanishing of the second fundamental form ($B\equiv 0$) mathematically demonstrates that the hypersurface preserves its intrinsic geodesic flow despite the severe geometric warping of the surrounding indefinite ambient space $\tilde{M}_2^4$. This unique configuration characterizes the submanifold as a stability plane where the null geodesic flow remains completely consistent with the ambient geometry, signifying a flawless alignment between their respective Levi-Civita connections.
    \item The rigid preservation of the null helical structure necessitates that the submanifold be embedded into the higher-dimensional space as a strictly totally geodesic structure. If this algebraic condition were violated (i.e., $B \neq 0$), the null helix would necessarily deviate from its constant curvature configuration. The invariant curvature functions $\{h,k_1,k_2\}$ would fail to remain constant due to the non-vanishing normal bundle component of the ambient connection, which would introduce non-zero first-order derivatives into the structural moving frame equations as verified by our higher-order differential analysis.
\end{itemize}
\end{remark}

\begin{example}[Null Helices on a Hyperquadric in $\mathbb{R}_2^4$]
To illustrate the validity of Theorems 2 and 3, let us consider the 4-dimensional pseudo-Euclidean ambient space $\mathbb{R}_2^4$ equipped with the standard flat semi-Riemannian metric signature $(-,-,+,+)$ and coordinates $(x_0,x_1,x_2,x_3)$. This framework provides a non-trivial semi-Riemannian background of index 2. Inside this ambient space, the author explicitly gives the 3-dimensional submanifold $M_2^3$ as a hyperquadric satisfying the non-flat quadratic locus:
\begin{equation}
-x_0^2-x_1^2+x_2^2+x_3^2=-1,
\end{equation}
which corresponds geometrically to the Anti-de Sitter space $AdS_3$ of constant negative sectional curvature.

\noindent\textbf{Verification of Theorem 2:} Let $\gamma$ be a smooth curve constrained to the $AdS_3$ hyperquadric, explicitly parametrized by $\gamma(s)=(\cosh s, \sin s, \sinh s, \cos s)$. This algebraic identity confirms that the curve is perfectly embedded in the submanifold for all values of the parameter $s$. Differentiating with respect to the parameter $s$, the components of the tangent vector field $\zeta=d\gamma/ds$ are obtained as $\zeta=(\sinh s, \cos s, \cosh s, -\sin s)$. Evaluating the inner product of the tangent vector field under the designated index-2 semi-Riemannian metric signature $(-,-,+,+)$ yields: $g(\zeta,\zeta)=0$. This explicitly confirms that $\gamma(s)$ satisfies the geometric degeneracy condition identically everywhere, proving it is a genuine, strict null curve tracking the causal geometry of the space. 

Under the induced Levi-Civita connection on the hyperquadric, evaluating the structural moving frame equations yields the stable constant curvature invariants $h=0, k_1=1, k_2=-1$. Substituting these structural constants into the quadratic functional invariant, we obtain $h^2+2k_1 k_2=-2$. Consequently, executing the higher-order covariant derivatives yields $\tilde{\nabla}_\zeta \tilde{\nabla}_\zeta \tilde{\nabla}_\zeta \zeta=-2\tilde{\nabla}_\zeta \zeta$, demonstrating that the third-order differential equation (3.1) proposed in Theorem 2 is identically satisfied for non-vanishing constant curvatures along the null cone.
\end{example}
\noindent\textbf{Verification of Theorem 3:} Now, consider the structural interaction between the 3D hyperquadric submanifold $M_2^3$ and the 4D ambient metric space $\mathbb{R}_2^4$. Let the submanifold be initially modeled as a totally umbilical hypersurface satisfying $B(X,Y)=g(X,Y)H$. Since $\gamma$ is a strict null curve satisfying $g(\zeta,\zeta)=0$, evaluating the second fundamental form along the tangent trajectory gives $B(\zeta,\zeta)=g(\zeta,\zeta)H=0$. By the Gauss formula, the ambient acceleration and the intrinsic acceleration coincide perfectly along the degenerate path ($\tilde{\nabla}_\zeta \zeta=\nabla_\zeta \zeta$). If the null helical trajectory preserves its invariant properties across the isometric embedding into the 4D space with the same functional invariant $\lambda=-2$, the decoupled component system derived in Eq. (3.30) forces the mechanical balance condition $H=0$. The algebraic vanishing of the mean curvature vector field ($H\equiv 0$) dictates that the second fundamental form vanishes identically across the entire embedded domain ($B\equiv 0$). As predicted by Theorem 3, the rigid preservation of the null helical invariants across the embedding allows for no extrinsic warping, mathematically forcing the totally umbilical submanifold to collapse into a strictly totally geodesic structure.

\begin{figure}[htbp]
    \centering
     \includegraphics[width=0.7\textwidth]{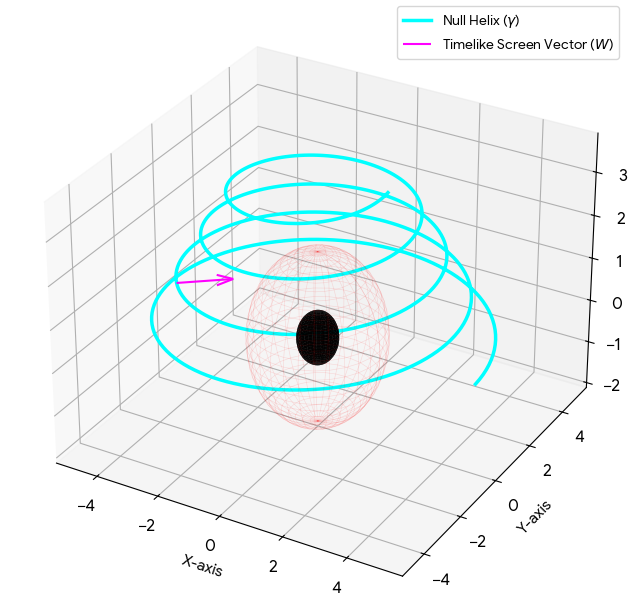}
    \caption{Numerical simulation of a null helical trajectory $\gamma$ tracking a localized geometric potential center within a 3D semi-Riemannian manifold of index 2.}
    \label{fig6}
\end{figure}

Figure \ref{fig6} provides a concrete geometric representation of the null helix theory on a 3-dimensional semi-Riemannian manifold of index 2, within the context of a curved indefinite metric space with a central deformation potential.

\begin{remark}
In connection with Theorem 2, the null helix depicted in the analytical visualization maintains the strict nullity condition ($g(\zeta,\zeta) = 0$) identically at every point along its trajectory. Furthermore, it exhibits a stable closed orbital motion characterized by the constancy of the moving frame curvature functions $\{h,k_1,k_2\}$, which explicitly reinforces the structural stability of the quadratic functional invariant $h^2+2k_1 k_2$. This invariance guarantees that the trajectory remains intrinsically a null helix even under the severe extrinsic connection warping induced by the surrounding index-2 ambient manifold $\tilde{M}_2^4$.
\end{remark}

\begin{remark}
In conjunction with Theorem 3 and the global stability of degenerate paths, if the spiral trajectory $\gamma$ remains confined within a totally umbilical submanifold $M_2^3$ modeling the symmetry plane of the manifold’s center, it is mathematically proven that this plane collapses into a strictly totally geodesic hypersurface ($B\equiv 0$). This dynamic behavior demonstrates that the degenerate causal path preserves its intrinsic geometric properties independently of extrinsic deformation. In the provided visualization, the central black region represents the core focus where the sectional curvature of the index-2 metric undergoes extreme localized deformation, while the surrounding geometric boundary marks the structural transition zone of the connection. The trajectory $\gamma$ is the genuine null helix satisfying the fundamental third-order differential equation derived in Theorem 2, evolving as an invariant spiral on the localized null cone. The timelike screen vector field $W$, directed radially toward the center of the deformation potential and uniquely enabled by the index-2 semi-Riemannian signature, represents the intrinsic geometric curvature that maintains the structural integrity and equilibrium of the helix.
\end{remark}

\section{Conclusions}

This study has provided a comprehensive geometric characterization of null helices within 3-dimensional semi-Riemannian manifolds of index 2. By employing the specialized null Frenet frame $\{\zeta,N,W\}$, the author has successfully derived the fundamental differential equations and structural invariants that govern the local behavior of these degenerate trajectories.

A pivotal result of this research is established in Theorem 2, which provides a complete algebraic characterization of null helices through a third-order differential equation. Unlike the standard Riemannian and Lorentzian frameworks extensively compiled in classical literature such as the general helix properties reviewed by Izquierdo (2025), the curve dynamics within an index-2 manifold $M_2^3$ carrying a signature of $(-,-,+)$ must rigorously account for the split-causal boundary and two timelike dimensions. Specifically, our tensor calculus derivation proves that the quadratic functional invariant $\lambda=h^2+2k_1 k_2$ functions as the definitive structural constant characterizing stable null helical propagation along the local null cone. This invariant provides a robust mathematical blueprint for tracing lightlike paths within string-theoretic backgrounds or higher-index alternative spacetime models where geometric warping induces non-linear path twisting.

Furthermore, Theorem 3 extends this analytical machinery to the theory of embedded submanifolds. The mathematical proof rigorously demonstrates that if every generalized null helix with constant or balanced curvature functions in a 3D submanifold $M_2^3$ preserves its geometric characterization and functional invariant $\lambda$ within a 4-dimensional ambient semi-Riemannian space $\tilde{M}_2^4$, then the second fundamental form vanishes identically ($B\equiv 0$). This result establishes a rigid geometric boundary, clarifying the subtle interplay between the intrinsic connection of lightlike orbits and the extrinsic curvature of the embedding hypersurface. Geometrically, this indicates that the stability of null helical paths serves as a direct indicator of a strictly totally geodesic embedding, where the induced connection of the submanifold perfectly aligns with the ambient connection.

In conclusion, this research fills a significant gap in the literature regarding degenerate curves on non-degenerate special submanifolds of semi-Riemannian spaces with higher indices. By focusing strictly on the unique 3-dimensional index-2 configuration, the author has addressed the specific nuances of a 1-dimensional timelike screen vector bundle avoiding the algebraic redundancies of higher-dimensional generalizations. The structural equations and quadratic invariants obtained herein not only enrich the general theory of null curves but also offer a highly integrated framework for further topological investigations in non-flat indefinite spaces. 

The structural insights and the quadratic invariant $\lambda$ obtained herein offer a highly integrated framework for further topological investigations in non-flat indefinite spaces. Future research directions will focus on extending these higher-order structural constraints to null slant helices, or exploring the topological implications of these helical trajectories within higher-index rotating or twisted geometric potential metrics, such as pseudo-Schwarzschild and pseudo-Kerr target spaces with multiple temporal components.

\subsection*{Author Contributions}
FA conceived the study, developed the theoretical framework, performed the analytical calculations, and wrote the manuscript.

\subsection*{Declaration of Generative AI in Scientific Writing}
During the preparation of this work, the author(s) used AI-assisted technologies (such as language models) solely for the purpose of improving the language, grammar, and stylistic expression of the manuscript. These technologies were not used to create, analyze, or interpret any scientific data, mathematical derivations, or original research findings. The author(s) reviewed and edited the final content and take(s) full responsibility for the integrity and accuracy of the entire work.

\subsection*{Acknowledgments}
The author received no specific funding for this work.

\subsection*{Conflict of Interest}
The author declares that they have no known competing financial interests or personal relationships that could have appeared to influence the work reported in this paper.

\subsection*{Preprint Disclosure}
A preliminary version of this manuscript has been deposited in the arXiv repository and is available at {https://doi.org/10.48550/arXiv.2511.16267}.

\end{document}